\documentclass[12pt]{article}

\usepackage{graphicx}
\usepackage{float}

\usepackage{caption}
\usepackage{subcaption}

\newtheorem{theorem}{Theorem}
\newtheorem{lemma}{Lemma}
\newtheorem{corollary}{Corollary}

\date{}

\include{abstract}

\begin{document}

\begin{titlepage}
\vspace{1.5in}
\begin{center}
\begin{Large}
{\bf {Random Knockout Tournaments  }}

\end{Large}

\vspace{0.2in}
\begin{large}
{\bf{Ilan Adler}} \\
Department of Industrial  Engineering and Operations Research\\
University of California, Berkeley\\
$\;$ \\
{\bf{Yang Cao}}\\
Department of Industrial and Systems Engineering\\
University of Southern California\\
$\;$ \\ 
{\bf{Richard Karp}} \\
Computer Science Department\\
University of California, Berkeley\\
$\;$\\
{\bf{Erol Pekoz}}\\
School of Business\\
Boston University\\
$\;$\\
{\bf{Sheldon M. Ross}}\\
Department of Industrial and Systems Engineering
$\;$\\
University of Southern California
$$\;$$
adler@ieor.berkeley.edu\\
cao573@usc.edu\\
karp@cs.berkeley.edu\\
pekoz@bu.edu \\
smross@usc.edu
$$\;$$

\end{large}

\vspace{0.2in}

\end{center}

\end{titlepage}

\begin{abstract} We consider a random knockout tournament among players $1, \ldots, n$, in which each match involves two players. The match format is specified by the number of matches played in each round, where the constitution of the matches in a round is  random. Supposing that there are numbers $v_1, \ldots, v_n$ such that  a match between $i$ and $j$ will be won by $i$ with probability $\frac{v_i}{v_i+v_j}$, we obtain a lower bound on the tournament win probability for the best player, as well as upper and lower bounds for all the players. We also obtain additional bounds by considering the best and worst formats for player $1$ in the special case $v_1 > v_2 = v_3 = \cdots = v_n.$ 

\end{abstract}

\section{Introduction}
We consider a tournament among players $1, \ldots, n$, in which each match involves two players.   The tournament is assumed to be of knockout type  in that the losers of matches are eliminated and do not move on to the next round, and the  tournament continues until all but one player is eliminated, with that player being declared the winner of the tournament. The match format is specified 
 by the set of positive integers  $r, m_1, \ldots, m_r$ with the interpretation that there are a total of $r$ rounds, with round $i$ consisting of $m_i$ matches, $\sum_{i=1}^r m_i = n-1.$ Because  $\sum_{j=1}^{i-1} m_j$ players  have been eliminated by the end of round $i-1$,  we must have that 
$m_{i} \leq (n - \sum_{j=1}^{i-1} m_j)/2.$

We  suppose that  the constitution of the matches in a round is totally random. That is, for instance, the $2m_1$ players that play in round $1$ are randomly chosen from all $n$ players and then randomly arranged into  $m_1$ match pairs. The winners of these $m_1$ matches, along with the $n - 2m_1$ players that did not play a match in round $1$ then move to round $2$, and so on.

We suppose that the players have respective values $v_1, \ldots, v_n$,  and that a match involving players $i$ and $j$ is won by player $i$ with probability $v_i/(v_i+v_j).$ We let 
$P_i $ be the probability that player $i$ wins the tournament, $i=1, \ldots, n.$  In Section 3 we derive a lower bound on the probability that the strongest player (e.g., the one with the largest $v$) wins the tournament, and an upper bound on the probability that the weakest player wins the tournament.  In Section 4 we derive upper and lower  bounds on the win probabilities 
$P_i, i \geq 1,$ and also show that if $v_1 \geq v_2 \geq \cdots \geq v_n$ then $P_1 \geq P_2 \geq \cdots \geq P_n.$   In Section 5, we consider the special case where  $v_1 > v_2 = \ldots = v_n$ and show that when $n = 2^s+ k, 0 \leq k < 2^s,$ the best format for the strongest player is the so-called {\it{balanced format}} that has $k$ matches in the first round and then has all remaining players competing in each subsequent round. We also show that whenever the number of  remaining players, say t, is even there is an optimal (from the point of view of the best remaining player) format that  calls for  $t/2$ matches in the next round. We also show, for the sections special case,  that the worst format for the best player is to have exactly one  match each round.  Analogous results for the worst player are also shown. Although we do not have a proof we conjecture that,  among all possible formats, the balanced format maximizes and the one match per round minimizes the best players probability of winning the tournament even in the case of general $v_i.$ (We show by a counterexample that the format that calls for $t/2$ matches when an even number $t$ of players remain is not optimal for the best player in the general case.)

\section{Literature Review}

\begin{figure}[H]
\centering
\begin{subfigure}{.5\textwidth}
  \centering
  \includegraphics[width=6cm]{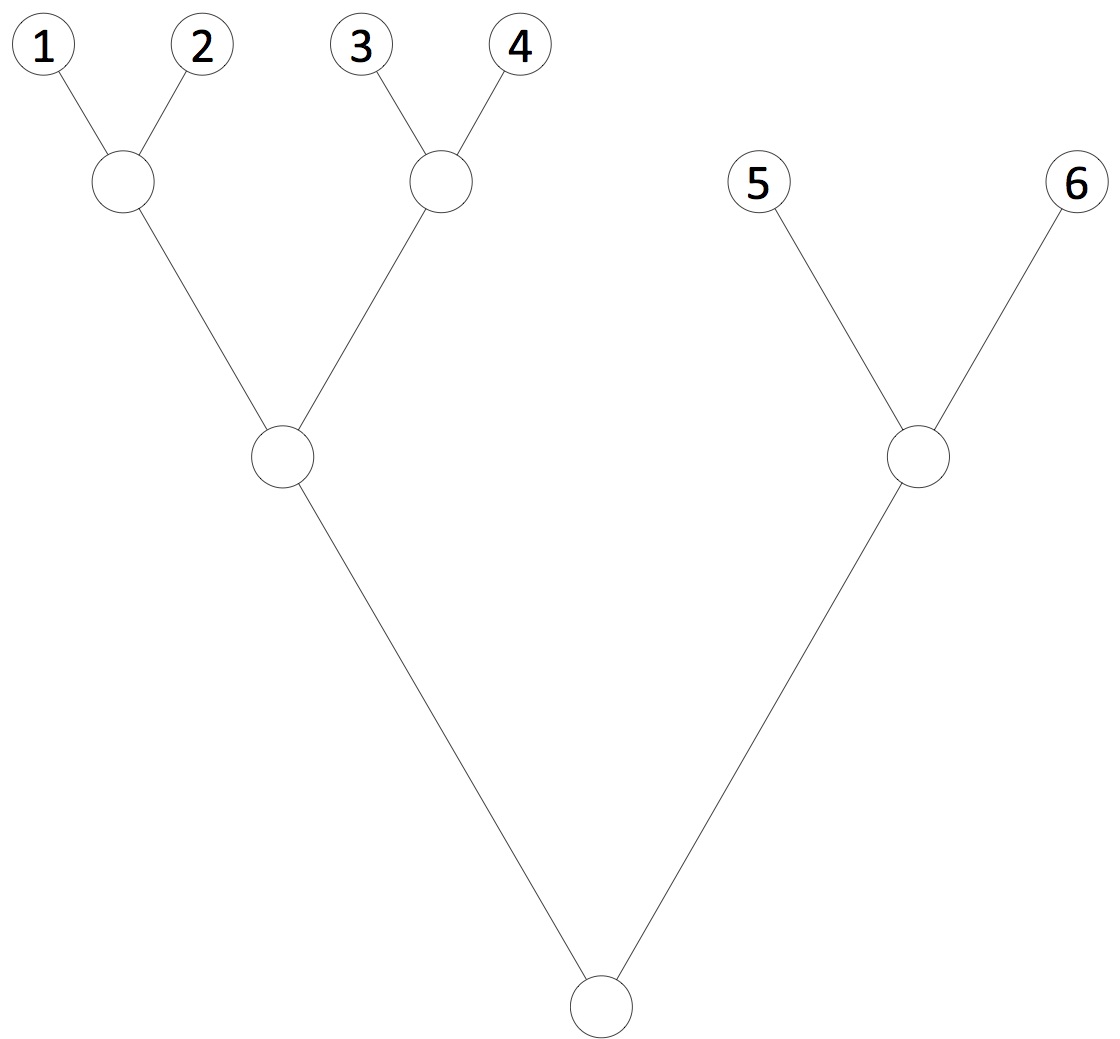}
  \caption{First balanced structure}
  \label{fig_bal_str1}
\end{subfigure}%
\begin{subfigure}{.5\textwidth}
  \centering
  \includegraphics[width=6cm]{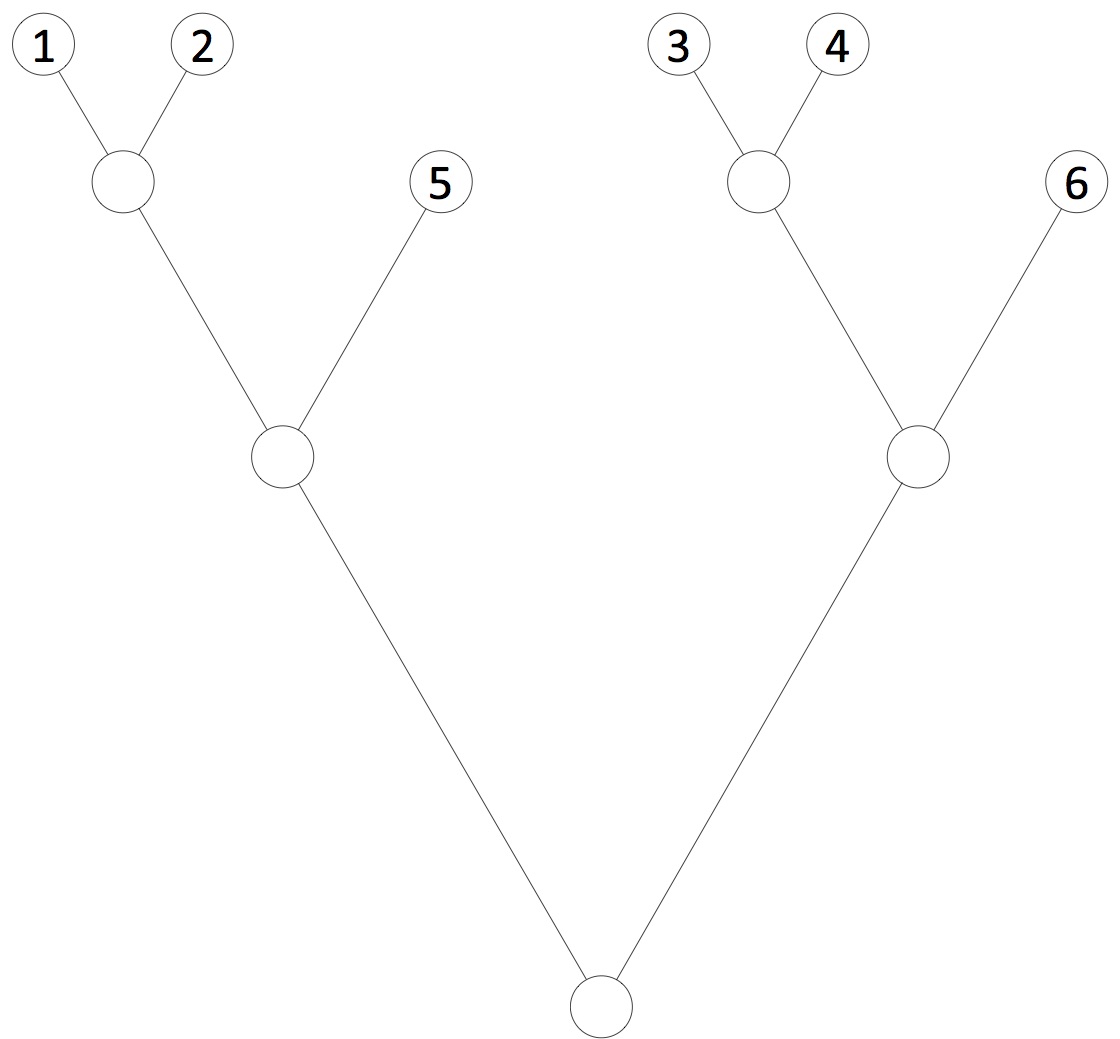}
  \caption{Second balanced structure}
  \label{fig_bal_str2}
\end{subfigure}
\caption{Two balanced structures with 6 players}
\label{fig_bal_str}
\end{figure}

Most papers in the literature on random knockout tournaments consider  more structured formats  than the ones we are considering  which suppose that the number of matches in each round is fixed and that  the game participants in each round are randomly chosen from those that remain. In these papers the structure of the tournament is fixed and players are randomly assigned to positions of the structure. For example, a structure with 6 players having a balanced format (meaning 2 matches in round 1, 2 matches in round 2 and 1 match in round 3) might in round 1 have a match between those in positions 1 and 2 and one between those in positions 3 and 4, with the winners then playing each other  in round 2  and 5 playing 6 in round 2.  Another possibility would be to have the winners of round 1  each play  in round 2 one of those that did not play in round 1. Figures \ref{fig_bal_str1} and \ref{fig_bal_str2} indicate these two structures.

An example of a  non-balanced structure with 6 players is one where there is a single  match in each round, with the winner of a match  playing in the following round
someone who has not yet played.

Maurer \cite{Maurer} proved  for  random structured formats, when $v_1>v = v_2 = v_3 = ... = v_n$,  that the win probability of player $1$ is maximized under the balanced structures. Because, as noted by Maurer, a random format can be expressed as a mixture of random structured formats, this also establishes the result for random formats. In section 5 we give another proof of this result, which is both quite elementary and also shows,   under the same condition,  that when there are $n=2^s+k$ remaining players the optimal random format can schedule either $k$  or $n/2$ games if $n$ is even, and  $k$ games if $n$ is odd. We also prove, when $v_1>v = v_2 = v_3 = ... = v_n$, that the win probability of  player $1$ is uniquely minimized by the random format that has exactly one match in each round.

Chung and Hwang \cite{ChungHwang} proved that for any structured format, if $v_i$ is decreasing in $i$ then so is the winning probability, which also proves the result for any random unstructured format. In Corollary 1 of Section 4 we give a proof of this result that  also proves, when  $v_i$ is decreasing in $i,$  that the probability player $i$ reaches round $s$ is decreasing in $i$ for all $s$.

There are other papers studying knockout tournaments, for example, David \cite{David}, Glenn \cite{Glenn}, Hennessy and Glickman \cite{HennessyGlickman}, Horen and Riezman \cite{HorenRiezman}, Hwang \cite{Hwang}, Marchand \cite{Marchand}. The paper  Ross and Ghamami \cite{RossGhamami} discusses efficient ways of simulating win probabilities for such tournaments.

\section{Lower Bound on the Strongest Player's Win Probability}
In this section we establish a lower bound on the probability that the strongest player wins the tournament, and an upper bound on the probability that the weakest player wins.

$$\;$$
\begin{theorem}
 If $v_i $ is decreasing in $i$, then 
$$P_1 \geq \frac{v_1}{\sum_{j=1}^n v_j}, \quad P_n \leq    \frac{v_n}{\sum_{j=1}^n v_j} .$$
\end{theorem}
$$\;$$
Preliminary to proving the theorem, we define the {\it{weight}} of a player as follows.  Say that a player is alive if that player has not been eliminated. 
If $S$ is the current set of alive players then, for $i \in S$,  the weight of player $i$ is defined to equal $v_i/\sum_{j \in S} v_j.$ If $i$ is no longer alive,  its weight is defined to be $0.$ 
$$\;$$
The following Lemma will be needed to prove the theorem. \\
\begin{lemma}
Suppose at the start of a round that players $A$ and $B$ are alive and   have respective  weights $x$ and $y$, where $x \geq y.$ If $X$ and $Y$ are random variables denoting the respective weights of $A$ and $B$  after the round, then $ y E[X] \geq x E[Y].   $   
\end{lemma} 
{\bf{Proof:}} The result is immediate if $y=0$ so suppose that $x \geq y > 0.$  Define  a random variable  $Z$ whose value is determined by the results of the following round, as follows: 
\begin{itemize}
\item $Z=0$ if $A$ and $B$ both lose
\item $Z=1 $ if $A$ and $B$ both advance
\item $Z=2 $ if $A$ and $B$ play each other
\item $Z= (i, j)$ if one of the players $A$ or $B$ defeats $i$ whereas the other loses to  $j$
\item  $Z= (0, j)$ if one of the players $A$ or $B$  advances by not being selected to play whereas the other loses to $j$
\end{itemize}
By considering the possible values of $Z$, we now show that  $y E[X|Z] \geq x E[Y|Z] .$ 
\begin{enumerate}
\item $\;E[X|Z=0] = E[Y|Z=0] = 0 $
\item  Conditional on $Z = 1$, let $R$ represent the sum of the weights of those aside from $A$ and $B$ that advance. Then, 
$$E[X|Z=1] = E[ \frac{x}{x+y+R}] , \quad E[Y|Z=1] = E[ \frac{y}{x+y+R}] $$
showing that  $\;yE[X|Z=1] = x E[Y|Z=1].$ 

\item  Conditional on $Z = 2$, let $R$ represent the sum of the weights of those aside from $A$ and $B$ that advance (and note that $R$ is independent of who wins the game between $A$ and $B$). Then, 
$$E[X|Z=2] = \frac{x}{x+y} E[\frac{x}{x+R}] , \quad E[Y|Z=2] = \frac{y}{x+y} E[\frac{y}{y+R}] $$
Hence, 
$$ y  E[X|Z=2] - x  E[Y|Z=2] = \frac{xy}{x+y}   E[\frac{x}{x+R} -  \frac{y}{y+R}  ]  \geq 0$$
where the inequality follows because $x \geq y$ implies that  $\frac{x}{x+r} \geq  \frac{y}{y+r} $ for any $r \geq 0.$

\item 
With $w$ and $v$ representing, respectively, the weights of $i$ and $j$ before the round, we have 
$$P(A \; \mbox{advances} | Z = (i,j) ) = c (\frac{x}{x+w} )   (\frac{v}{v+y} )   $$
where  $\; 1/c =  (\frac{x}{x+w} )   (\frac{v}{v+y} ) +  (\frac{y}{y+w} )   (\frac{v}{v+x} ) .$ Let $R$ represent the sum of the weights of those, aside from $A, B$ and $j$, that advance, and note that $R$ is independent of which of $A$ or $B$ loses. Then
$$ E[X|Z= (i,j)]  = c(\frac{x}{x+w} )   (\frac{v}{v+y} ) E[ \frac{x}{R+x+ v} ] $$
and, similarly,

$$ E[Y|Z= (i,j)]  = c(\frac{y}{y+w} )   (\frac{v}{v+x} ) E[ \frac{y}{R+y+ v} ] $$
Hence, letting  $D = y E[X|Z= (i,j)]  - x  E[Y|Z= (i,j)] $, we see that 

$$ D  = c\,x\,y\,v \,E[  \frac{x}{ (x+w)(y+v)(R + x + v) } \, - \, \frac{ y}{ (y+w)(x+v) (R+y+v)} ]$$
As it is easy to check that for every $R \geq 0$ and $x \geq y$
$$x(y+w)(x+v) (R+y+v) \geq y (x+w)(y+v)(R + x + v) $$
it follows that 
 $y E[X|Z= (i,j)]  - x  E[Y|Z= (i,j)] \geq 0 $ for  $x \geq y.$ 
\item That  $y E[X|Z= (0,j)]  - x  E[Y|Z= (0,j)] \geq 0 $ follows from the preceding result by setting $w=0.$ 

\end{enumerate}
Hence, we have shown that  $\;y E[X|Z]  \geq x E[Y|Z],$ and the result follows upon taking expectations of both sides of this inequality. $\qquad \rule{2mm}{2mm}$
$$\;$$
{\bf{Proof of Theorem 1:}}  We give the proof that $P_1 \geq v_1/ \sum_{j=1}^n v_j.$ The proof that  $P_n \leq v_n/ \sum_{j=1}^n v_j$ is similar.  Let $W_j(k)$ be the weight of player $j$ after $k$ rounds have been played. Also, let $H_k$ be the history of all results through the first $k$ rounds. We claim that
$$E[W_1(k+1)| H_k] \geq W_1(k)$$
Because the claim is true when $W_1(k) = 0,$ assume that $W_1(k) > 0.$ Let $A_k$ denote the set of alive players after round $r$. Now, from Lemma 1, we have, for $j \in A_k,$ that
$$W_j(k) E[W_1(k+1)| H_k] \geq W_1(k) E[W_j(k+1)| H_k] $$
Hence, if $E[W_1(k+1)| H_k] < W_1(k)$, then for any $j \in A_r$
$$ W_j(k) > E[W_j(k+1)| H_k] $$
which is a contradiction since  $1 = \sum_{j \in A_k} W_j(k) =  \sum_{j \in A_k}   E[W_j(k+1)| H_k] .$
Hence, $E[W_1(k+1)| H_k] \geq W_1(k)$, and taking expectations of both sides gives 
$$ E[W_1(k+1)] \geq E[W_1(k)], \; r \geq 0.$$
If the tournament has $r$ rounds, then the preceding yields that $E[W_1(r)] \geq E[W_1(0)]$, which gives the result since $E[W_1(r)] = P_1$ whereas $W_1(0) = v_1/\sum_{j=1}^n v_j.
\qquad \rule{2mm}{2mm}$

$$\;$$
{\em{Remark:}} The preceding argument shows that $W_1(k),  k \geq 0$ is a submartingale, and that $W_n(k),  k \geq 0$ is a supermartingale.

$$\;$$
\section{Bounds on Win Probabilities}
Let 
$$p_i = \frac{1}{n-1} \sum_{j \neq i} \frac{v_i}{v_i + v_j}$$ be the probability that $i$ would win a match against a randomly chosen opponent. In this section we prove that $P_i$ is smaller than it would be if it were the case that $i$ would win each game it plays with probability $p_i.$ That is, we will prove the following.

\begin{theorem} If the tournament format is $(r, m_1, \ldots, m_r)$ then
$$P_i \leq \prod_{s=1}^r  (\frac{2m_s}{r_s} p_i + 1 -  \frac{2m_s}{r_s} )$$
where    $r_s = n - \sum_{j=1}^{s-1} m_j$ is the number of players that advance to round $s$.
\end{theorem}
To prove the preceding theorem we will need a couple of lemmas. Before giving these lemmas, we introduce the following notation. 
We let $R_{i,s}$ be the event that player $i$ reaches round $s, \, s = 1, \ldots, r$,  and let  $R_{i,r+1}$ be the event that $i$ wins the tournament.   If a player receives a bye in a round (that is, if it reaches that round but is not chosen to play a match) say that it plays player $0$ in that round.
Also, let $p_{ij} = \frac{v_i}{v_i+v_j},\;i\neq j,$ be the probability that $i$ beats $j$ in a game.

Lemma 2 is easily proven by a coupling argument.

\begin{lemma}  \label{StrDecLem}
  For all $s=1, \ldots, r+1, \, P(R_{i,s}R_{j,s})$ is  an increasing function of $v_j.$ (When $i=j$, this states that $P(R_{i,s})$ is an increasing function of $v_i.$)
\end{lemma}

\begin{lemma}  Assume $v_1 \geq v_2 \geq ... \geq v_n$. 
For fixed $i$, $P(R_{i,s}R_{j,s})$ decreases in $j$ for $j \neq i$.
\end{lemma}
{\bf{Proof:}}  Fix $i, j, k$, where $j < k$ and $j, k \neq i.$ The proof is by induction on $n$.  For $n=3$, the only format is to play one game each round. The result holds in this case because $ v_j \geq v_k$  implies that 
$$  3 P(R_{i,2}R_{j,2}) =\frac{v_i}{v_i + v_k} + \frac{ v_j}{v_j + v_k}  \geq \frac{v_i}{v_i + v_j} + \frac{ v_k}{v_j + v_k} = 3 P(R_{i,2}R_{k,2}) $$
Assume the results holds for up to $n-1$ players and for all  formats. We now consider the $n$ player case. Let $\pi$ be an arbitrary  format. Define a random vector $Z$ whose value is determined by the results of the first round of the tournament under $\pi$, as follows:

\begin{itemize}
\item $Z = (1,u)$ if $i$ plays against $j$ and $k$ plays against $u$, or $i$ plays against $k$ and $j$ plays against $u$;
\item $Z = (2,u)$ if $i$ plays against $u,$ and $j$ plays against $k$;
\item $Z = (3,u,v,w)$ if $i$ plays against $u$, $j$ beats $v,$ and $k$ beats $w$;
\item $Z = (4,u,v,w)$ if $i$ plays against $u$, $j$ loses to $v,$ and $k$ loses to $w$;
\item $Z = (5,u,v,w)$ if either (a): $i$ plays against $u$, $j$ beats $v,$ and $k$ loses to $w$, or (b):  $i$ plays against $u$, $k$ beats $v,$ and $j$ loses to $w$.
\end{itemize}

By considering the possible values of Z, we now show that $P(R_{i,s}R_{j,s}|Z) \geq P(R_{i,s}R_{k,s}|Z)$.

\begin{enumerate}
\item Let  $p_{zy} = \frac{v_z}{v_z+v_y}$.  Because  $p_{ik}\geq p_{ij}$ and $p_{ju}\geq p_{ku}$, it follows that
$$P(R_{i,2}R_{j,2}|Z = (1,u)) = p_{ik}p_{ju}/2 \geq  p_{ij}p_{ku}/2 = P(R_{i,2}R_{k,2}|Z = (1,u)).$$
Let $T$ denote the set of players other than $i, j, k,$  that reach round 2. Given that $i$ and $j$ have reached round 2, the probability of $i$ and $j$ reaching round $s>2$ equals to the probability that $i$ and $j$ reaches round $s-1$ in a new tournament which begins with players $i$, $j$ and $T$,  and follows the same format as $\pi$ after round 1. Let $P_{i,j,T}(s)$ denote the probability that $i$ and $j$ reaches round $s-1$ in the new tournament. Similarly, let $P_{i,k,T}(s-1)$ denote the probability that $i$ and $k$ reaches round $s-1$ in a tournament which begins with players $i$, $k$ and $T$, and follows the same format as $\pi$ after round 1. Then by Lemma \ref{StrDecLem}, $P_{i,j,T}(s-1)\geq P_{i,k,T}(s-1)$. Therefore for $s > 2$,
\begin{eqnarray*}
P(R_{i,s}R_{j,s}|Z = (1,u) ) &=& P(R_{i,2}R_{j,2}|Z = (1,u))E[P_{i,j,T}(s-1) ]  \\
&  \geq & P(R_{i,2}R_{k,2}|Z = (1,u))E[P_{i,k,T}(s-1) ]  \\
&=& P(R_{i,s}R_{k,s}|Z = (1,u)).
\end{eqnarray*}
where the preceding used that $T$ is independent of the event $R_{i,s}R_{j,s} $  and of the event $R_{i,s}R_{k,s} .$

\item Since $v_j\geq v_k$, we have that
$$P(R_{i,2}R_{j,2}|Z = (2,u)) = p_{iu}p_{jk} \geq  p_{iu}p_{kj} = P(R_{i,2}R_{k,2}|Z = (2,u)) $$
Let $T$ denote the set of players other than $i, j$ and $k$ that reaches round 2. With the same definition and argument as when $Z=(1,u)$, it can be shown that for $s>2$,
$$P_{i,j,T}(s)\geq P_{i,k,T}(s)$$
which implies that 
$$P(R_{i,s}R_{j,s}|Z = (2,u))  \geq P(R_{i,s}R_{k,s}|Z = (2,u)).$$

\item In this case,
$$P(R_{i,2}R_{j,2}|Z = (3,u,v,w)) = p_{iu} = P(R_{i,2}R_{k,2}|Z = (3,u,v,w)).$$
Now, given that $i, j$ and $k$ have reached round 2, the induction hypothesis implies that the probability of $i$ and $j$ reaching round $s$ is greater than or equal to the probability of $i$ and $k$ reaching round $s$ for $s>2$. Therefore for $s>2$,
$$P(R_{i,s}R_{j,s}|Z = (3,u,v,w)) \geq P(R_{i,s}R_{k,s}|Z = (3,u,v,w)).$$

\item In this case, $P(R_{i,s}R_{j,s}|Z = (4,u,v,w)) = P(R_{i,s}R_{k,s}|Z = (4,u,v,w)) = 0$.

\item In this case,
$$P(R_{i,2}R_{j,2}|Z = (5,u,v,w)) = c p_{iu} p_{jv}p_{wk}$$
$$P(R_{i,2}R_{k,2}|Z = (5,u,v,w)) = c p_{iu} p_{kv}p_{wj}$$
where $1/c = p_{jv}p_{wk}+p_{kv}p_{wj}$.

Since $p_{jv}\geq p_{kv}$ and $p_{wk}\geq p_{wj}$, we have that 
$$P(R_{i,2}R_{j,2}|Z = (5,u,v,w))\geq P(R_{i,2}R_{k,2}|Z = (5,u,v,w))$$
Let $T$ denote the set of players other than $i, j$ and $k$  that reach round 2, and note that $T$ has the same distribution whether (a) or (b) resulted. With the same definition and argument as when $Z=(1,u)$, it can be shown that for $s>2$,
$$P_{i,j,T}(s)\geq P_{i,k,T}(s)$$
which implies that 
$$P(R_{i,s}R_{j,s}|Z = (5,u,v,w))  \geq P(R_{i,s}R_{k,s}|Z = (5,u,v,w)).$$

\end{enumerate}
Hence, we have shown that $P(R_{i,s}R_{j,s}|Z) \geq P(R_{i,s}R_{k,s}|Z)$, and the result follows upon taking
expectations of both sides of this inequality.  $\qquad \rule{2mm}{2mm}$ \\
\vspace{.1mm}\\
{\bf{Proof of Theorem 2:}}  Assume that $v_j$ is decreasing in $j$. Now, 
given that $i$ reaches round $s$, the conditional probability that $j$ also reaches round $s$ is
$$P(R_{j,s}|R_{i,s}) = \frac{  P(   R_{j,s} R_{i,s}) }{  P(R_{i,s}) }$$
Hence, from Lemma  3 it follows that $P(R_{j,s}|R_{i,s}) $ is a decreasing function of $j, j \neq i$. Now, 
$$P_i(v_1, \ldots, v_n) = P(R_{i,2} \cdots R_{i,r+1} ) = \prod_{s=1}^r P(R_{i,s+1}| R_{i,s} )  $$
  Let $C_{i,s}$ be the event that $i$ competes in round $s$ (that is, it is the event that $i$ reaches round $s$ and then plays a match in that round).  With 
   $Q_j =   P( i \;\mbox{plays} \; j \; \mbox{in round}\;s| C_{i,s}), j \neq i,$ we have that

  \begin{eqnarray*}
  Q_j  & = &P( i \;\mbox{plays} \; j \; \mbox{in round}\;s| C_{i,s} R_{j,s}) P(R_{j,s}|C_{i,s}) \\
   & = & P( i \;\mbox{plays} \; j \; \mbox{in round}\;s| C_{i,s} R_{j,s})  P(R_{j,s}|R_{i,s})  \\
     & = & \frac{1}{r_s-1} P(R_{j,s}|R_{i,s}) 
       \end{eqnarray*}
       where  $r_s = n - \sum_{j=1}^{s-1} m_j.$  Hence, 
     $Q_j, j \neq i$ is a decreasing function of $j$. Letting $Y$ be a random variable such that $P(Y = j) = Q_j, j \neq i$, it thus follows that $Y$ is stochastically smaller than the random variable $X$ having $P(X=j) = 1/(n-1), j \neq i.$ 
Therefore,
       \begin{eqnarray*}
       P(R_{i,s+1}| R_{i,s} )  &=& 1 - \frac{2m_s}{r_s} +   \frac{2m_s}{r_s} P(R_{i,s+1}| C_{i,s} )  \\
        &=& 1 - \frac{2m_s}{r_s} +   \frac{2m_s}{r_s} \sum_{j \neq i} \frac{v_i}{v_i + v_j} Q_j \\
         &=& 1 - \frac{2m_s}{r_s} +   \frac{2m_s}{r_s}E[\frac{v_i}{v_i + v_Y} ] \\
          &\leq& 1 - \frac{2m_s}{r_s} +   \frac{2m_s}{r_s} E[\frac{v_i}{v_i + v_X} ] \\
          &=& 1 - \frac{2m_s}{r_s} +   \frac{2m_s}{r_s} \frac{1}{n-1} \sum_{j \neq i}  \frac{v_i}{v_i + v_j} 
          \end{eqnarray*}
        where the   preceding used that  $\frac{v_i}{v_i + v_j}$ is an increasing function of $j$ and $X$ is stochastically larger than $Y$ to conclude that $E[\frac{v_i}{v_i + v_Y} ] \leq E[\frac{v_i}{v_i + v_X} ].  \qquad \rule{2mm}{2mm}$

$$\;$$
{\em{Remark:}}
It follows from Lemma 2 that for all random formats, $P_i = P(R_{i, r+1})$ is an increasing function of $v_i$. On the other hand, it seems intuitive that $P_i$ is a decreasing function of $v_j$ for $j\neq i$. However,  while this is true for $n = 3$, it is not true for $n\geq 4$. The argument, when $n=3$, uses that

$$P_1 = \frac{1}{3}(\frac{v_2}{v_2+v_3}\frac{v_1}{v_1+v_2}+\frac{v_3}{v_2+v_3}\frac{v_1}{v_1+v_3})+ \frac{2}{3}\frac{v_1}{v_1+v_2}\frac{v_1}{v_1+v_3}.$$
This gives 
\begin{eqnarray*}
\frac{\partial P_1}{\partial v_2} &=&\frac{1}{3}v_1 \frac{\partial}{\partial v_2} \frac{v_2}{(v_2+v_3)(v_1+v_2)}+\frac{1}{3}\frac{v_1 v_3}{v_1+v_3} \frac{\partial}{\partial v_2} \frac{1}{v_2+v_3}+ \frac{2}{3}\frac{v_1^2}{v_1+v_3} \frac{\partial}{\partial v_2} \frac{1}{v_1+v_2}\\ 
&=& \frac{1}{3}\frac{v_1}{(v_1+v_2)^2 (v_2+v_3)^2(v_1+v_3)}\bigg( -v_1v_3^2 - 2v_2^2 v_3 -3v_1v_2^2 -6v_1 v_2 v_3 \bigg)\\
&\leq & 0.
\end{eqnarray*}

For $n=4$ a counter example can be constructed as follows. Consider the balanced format and let $P_i(v_1,v_2, v_3,v_4)$ denote the probability that $i$ wins the tournament when player $j$ has value $v_j, j=1,2,3,4$.  Conditioning on whether or not player 1 first plays against player 4, we have
$$P_1(2,1,1,x)= \frac{1}{3}\left(\frac{2}{2+x} \cdot \frac{2}{3}\right)+\frac{2}{3}\cdot \frac{2}{3}
\left(\frac{1}{1+x} \cdot \frac{2}{3}+\frac{x}{(1+x)}\frac{2}{(2+x)}\right)$$
and thus
$$P_1(2,1,1,\frac{1}{100})=\frac{31600}{60903}\approx .518858\ldots < .518861\ldots \approx \frac{7744}{14925}=P_1(2,1,1,\frac{1}{99}).  \qquad \rule{2mm}{2mm}$$


$$\;$$

Using Lemma 3, it is easy to show that 
 if $v_i$ is decreasing in $i$ then so is $P(R_{i,s})$.\\

\begin{corollary}
If $v_1 \geq ... \geq v_n$ then $P(R_{1,s})\geq P(R_{2,s}) \geq ... \geq P(R_{n,s}), \; s=1, \ldots, r+1$.
\end{corollary}
{\bf{Proof:}} Suppose  $i<j$. Because $r_s$ is the number of players that are still alive at the start of round $s$, it follows that given $i$ reaches round $s$, the expected number of others that also reach round $s$ is $r_s - 1.$ Consequently,
$$r_s - 1 =  \sum_{k \neq i} P(R_{k,s}|R_{i,s}) $$
Hence, 

\[
  (r_s-1)P(R_{i,s})=\sum_{k \neq i,j}P(R_{i,s}R_{k,s}) + P(R_{i,s}R_{j,s}),
\]
\[
  (r_s-1)P(R_{j,s})=\sum_{k \neq j,i}P(R_{j,s}R_{k,s}) + P(R_{j,s}R_{i,s}),
\]
which by Lemma 3 shows the desired result for  $s \leq r$. In addition, we have, 
\[
     P(R_{i,r+1})=\sum_{k \neq j,i}P(R_{i,r}R_{k,r}) \frac{v_i}{v_i+v_k}+ P(R_{j,r}R_{i,r})\frac{v_i}{v_i+v_j},
\]
\[
     P(R_{j,r+1})=\sum_{k \neq j,i}P(R_{j,r}R_{k,r}) \frac{v_j}{v_j+v_k}+ P(R_{i,r}R_{ijr})\frac{v_j}{v_j+v_i},
\]
which by Lemma 3, and by the assumption that $v_i \geq v_j$, completes the proof. $\qquad \rule{2mm}{2mm}$ \\
\vspace{.1mm}\\

We now give lower bounds on the win probabilities. 

\begin{theorem} 
Suppose $v_1 \geq v_2 \geq \ldots \geq v_n.$ Let  $X_1, \ldots, X_r$ be independent with $P(X_j = 1) =   \frac{2m_j}{r_j} = 1 - P(X_j = 0 ) , j = 1, \ldots,r,$ and let $N= \sum_{j=1}^r X_j$. Then
$$P_i \geq \sum_{g=1}^{i-1}   P(N = g)  \prod_{k=1}^{g} \frac{v_i}{v_i + v_k}  + \sum_{g=i}^{r}   P(N = g)   \prod_{k=1}^{i-1} \frac{v_i}{v_i + v_k}  \prod_{k=i+1}^{g+1} 
\frac{v_i}{v_i + v_k} $$
\end{theorem}
{\bf{Proof:}} 
Define the surrogate of $i$ as follows:
The initial surrogate of $i$ is $i$ itself; 
anyone who beats  a surrogate of $i$  becomes the current surrogate of $i$. 
Note that at any time there is exactly one player who is the current surrogate of $i$. If we let 
$X_j$ be the indicator of whether the surrogate of $i$ plays a match in round $j$, it follows that $X_1, \ldots, X_r$ are independent with $P(X_j = 1) =   \frac{2m_j}{r_j} = 1 - P(X_j = 0 ) .$ Also,  let $N= \sum_{j=1}^r X_j$ be  be the number of games played by  surrogates of $i$ (while they are the current surrogate),  and let $O = \{J_1, \ldots J_N\}$ be their set of  opponents in these games.  Because $\;P(i \; \mbox{wins the tournament}|N, O) = \prod_{k=1}^N \frac{v_i}{v_i + v_{J_k}},$ it follows that 
\begin{eqnarray*}
P_i &=& E[   \prod_{k=1}^N \frac{v_i}{v_i + v_{J_k}} ] \\
& \geq & E[   \prod_{k=1}^{\min(N, i-1)} \frac{v_i}{v_i + v_k}   \prod_{k=i+1}^{N+1} \frac{v_i}{v_i + v_k} ]
\end{eqnarray*}
where  $ \prod_{k=i+1}^{N+1} \frac{v_i}{v_i + v_k} $ is equal to $1$ when $N < i.$ Hence, 
$$P_i \geq \sum_{g=1}^{i-1}   P(N = g)  \prod_{k=1}^{g} \frac{v_i}{v_i + v_k}  + \sum_{g=i}^{r}   P(N = g)   \prod_{k=1}^{i-1} \frac{v_i}{v_i + v_k}  \prod_{k=i+1}^{g+1} 
\frac{v_i}{v_i + v_k}   \qquad \rule{2mm}{2mm}$$
{\em{Remarks:}}
\begin{enumerate}
\item 
The probability mass function of $N$, which does not depend on $i$, is easily obtained by solving recursive equations (see  Example 3.24 of Ross \cite{sr}).
 In the special case where $n = 2^s + k,\; 0 \leq k < 2^s,$ and where there are $k$ matches in round $1$ and afterwards all remaining players have matches in each subsequent round, we have
$$P(N= s+1) = \frac{2k}{n} = 1 - P(N=s).$$
\item A weaker bound than the one provided in   is obtained by noting that, for $i>1$,  $\frac{v_i}{v_i + v_{J_k}} \geq \frac{v_i}{v_i+ v_1} $ and so for $i>1$
\begin{eqnarray*}
P_i &\geq & E[ (\frac{v_i}{v_i+ v_1} )^N] \\
&=&  E[ (\frac{v_i}{v_i+ v_1} )^{X_1 + \ldots X_r}] \\
&=& \prod_{j=1}^r E[ (\frac{v_i}{v_i+ v_1} )^{X_j}] \\
&=& \prod_{j=1}^r (  \frac{2m_j}{r_j}  \frac{v_i}{v_i+ v_1} + 1- \frac{2m_j}{r_j}  )
\end{eqnarray*}
\end{enumerate}

\section{Special Case Best and Worst Formats for the Strongest  and Weakest Players} 

\begin{theorem}
Suppose $v_1 > v = v_2 = v_3 = \ldots = v_n.$ If  $n = 2^s + k, \, 0 \leq k < 2^s,$ then the balanced format that schedules $k$ matches  in round $1$ and then has all remaining players competing in each subsequent round leads to the maximal possible value of $P_1.$ 
\end{theorem}
{\bf{Proof:}} 
Let $p = \frac{v_1}{v_1+v} > .5$ be the probability that player $1$ wins in a match against another player. 
The proof is by induction. As there is nothing to prove when $n=2$, let $n = 2^s + k,  \, 1 \leq k \leq 2^s,$ and suppose that the result is true for all smaller values of $n$. Consider a format that calls for $j$  matches in the first round, where $j < k.$ Then, by the induction hypothesis,  the format of this type that is best for player $1$ will call for $k-j$ matches in the second round. The probability that player $1$ is among the final $2^s$ players under these conditions is 
$$f_1(p) \equiv [\frac{2j}{n} \,p + 1 - \frac{2j}{n}  ]   [\frac{2(k-j) } {n-j} \,p + 1 - \frac{2(k-j)}{n-j} ]$$
However, the probability that player $1$ is among the final $2^s$ players if the first round has $k$ matches  is 
$$f_2(p) \equiv \frac{2k}{n} \,p + 1 - \frac{2k}{n} $$
Now, $f_1(p) - f_2(p) =0$ when $p$ is either $1/2$ or $1$. Because $f_1(p) - f_2(p)  $ is easily seen to be strictly convex this implies that $f_1(p) - f_2(p)  < 0$ when $1/2 < p < 1.$ Thus, in searching for the best format for player $1$ we need not consider any format that calls for less than $k$  matches in round $1$.

Now consider any format that calls for $k+i$  matches in round $1$ where $i > 0$. By the induction hypothesis it will call for $ 2^{s-1} - i $ matches in round $2$. This will result in  player $1$ being among the final $2^{s-1}$ players with probability  
$$  g_1(p)  \equiv [\frac{2k+2i} {n} \,p + 1 - \frac{2 k + 2i}{n}  ]   [\frac{ 2^s - 2i } {2^s - i} \,p + 1 - \frac{ 2^s - 2 i } {2^s - i} ]$$
However, the format that calls for $k$  matches in round $1$ and $2^{s-1}$  matches in round $2$ leads to player $1$ being  among the final $2^{s-1}$ players with probability  
$$g_2(p)    \equiv (\frac{2k}{n} \,p + 1 - \frac{2k}{n})  p$$
Now, $g_1(p) - g_2(p) =0$ when $p$ is either $1/2$ or $1$. 
 Because 
$$g_1''(p)=     2 \,\frac{2k+2i} {n}     \,\frac{ 2^s - 2i } {2^s - i}   , \quad   g_2''(p)  = 2\, \frac{2k}{n} ,$$          
it follows that  $ g_1''(p) \geq g_2''(p) $ is equivalent to 
$$(k+i) (2^s - 2i) \geq k (2^s-i) ,$$
which  is easily seen to be equivalent to 
$$2^s \geq k+2i,$$
which holds because having $k+i$ two person matches implies that $2^s + k \geq 2(i+k).$  Hence, $g_1(p) - g_2(p)$ is convex, which shows that there is an optimal format that initially has $k$ matches in round 1. The induction hypothesis then proves the result. $\qquad \rule{2mm}{2mm}$ 
$$\;$$
{\em{Remarks}} \\
1. Because  $g_1(p) - g_2(p)$ is strictly convex unless $2^s+ k = 2(i+k)$, it follows that the format specified in Theorem 2 is uniquely optimal (in the sense of maximizing $P_1$) when $n$ is odd, whereas when $n$ is even there is also an optimal format that schedules $n/2$ matches in round 1.\\
\vspace{.1mm}\\
2. In the case of general $v_i$, the format that calls for $n/2$ matches when an even number $n$ of players remain is not optimal for the best player. For a counter example, consider a knockout tournament with 6 players having values $v_1 = 6, v_2 = 4, v_3 = 3, v_4 = v_5 = v_6 = 1$. Under the balanced format, the probability of player 1 winning the tournament is $P_1 = 0.4422$. Under a format that plays 3 games in round 1 and then one game in each round, the probability of player 1 winning the tournament is $P_1' = 0.4412 < P_1$.
$$\;$$
\begin{theorem}
 Suppose $v_1 > v = v_2 = v_3 = \ldots = v_n.$ The unique format that minimizes $P_1$ is the one that has exactly one match in each round.
\end{theorem}
{\bf{Proof:}} The proof is by induction. Suppose it is true for all tournaments with fewer than $n$ players and now suppose there are $n$ players. Consider any format that calls for $s$ matches in the first round, where $s>1$. The probability that player $1$ is still alive when there are only $n-s$ alive players is
$$f_1(p) = \frac{2s}{n} p + 1 - \frac{2s}{n} $$
On the other hand, the probability that player $1$ is still alive when there are only $n-s$ alive players when the format is one match per round is
$$f_2(p) = \prod_{j=0}^{s-1} (\frac{2}{n-j} p + 1 - \frac{2}{n-j})  $$
Because $f_2(p)$ is a polynomial whose  coefficients are all positive it follows that $f_2''(p) > 0.$ As $f_1''(p) = 0$ this implies that $f_1(p) - f_2(p)   $    is strictly concave, which since 
$f_1(p) - f_2(p)  =0 $  when $p=1/2$ or $p=1$, enables us to conclude that $f_1(p) - f_2(p)  > 0$ when $1/2 < p < 1.$ Hence,  by the induction hypothesis, the format that calls for 
one match in each round results in a win probability for player 1 that is strictly smaller than what it is under any format that calls for $s > 1$ matches in round 1.  $\qquad \rule{2mm}{2mm}$ 
$$\;$$

The following theorem gives the analogous results for the weakest player. The proofs are similar and thus omitted.
\begin{theorem}
Suppose $v_1  = v_2 = v_3 = \ldots = v_{n-1} > v_n.$ The format that results in the highest win probability for player $n$ is the one that has one match in each round. 
 If $n = 2^s + k, \, 1 \leq k \leq 2^s,$ then the format that minimizes $P_n$ is the one that schedules $k$ matches  in round $1$ and then has all remaining players competing in each subsequent round. 
 \end{theorem}

We now use the preceding results to obtain universal (that is, they hold for all formats) upper bounds on $P_i.$  Recall that  $\, p_i = \frac{1}{n-1} \sum_{j \neq i} \frac{v_i}{v_i + v_j}.$\\

\begin{lemma} Suppose $n = 
 2^s+ k, \,0 \leq k < 2^s$.
 {\it
 \begin{description}
  \item[(i)]
   If $p > \frac{1}{2}$, then   $\prod_{j=1}^r  \left(\frac{2m_j}{r_j} p + 1 -  \frac{2m_j}{r_j} \right)$ is maximized by the balanced format, and
   its value is  $ \left(\frac{2k}{n}p + 1 - \frac{2k}{n} \right) p^s$.
   \item[(ii)]
   If $p = \frac{1}{2}$, then   $\prod_{j=1}^r   \left(\frac{2m_j}{r_j} p + 1 -  \frac{2m_j}{r_j}  \right)=\frac{1}{n}$ for all eligible formats.
   \item[(iii)]
   If $p < \frac{1}{2}$,  $\prod_{j=1}^r   \left(\frac{2m_j}{r_j} p + 1 -  \frac{2m_j}{r_j}  \right)$ is maximized by the format where
   at each round there is exactly one match, and
   its value is  $ \prod_{j=1}^{n-1} \left(\frac{2}{n-j+1}p + 1 - \frac{2}{n-j+1} \right)$.
   \end{description}}
   \end{lemma}
   {\bf Proof.} \hspace{1mm} Note that for a tournament where the players have respective values $(v_1,v \ldots, v)$,
   $P_1=\prod_{j=1}^r  \left(\frac{2m_j}{r_j} p_1 + 1 -  \frac{2m_j}{r_j} \right)$. Now, (i) and (iii) are direct corollaries of Theorems 2 and 4 respectively, while (ii) is true because
   if $v_1=v$ then all players have equal probability of winning the tournament regardless of the format.  $\qquad \rule{2mm}{2mm}$ \\
   \begin{theorem} Suppose $n = 
 2^s+ k, \,0 \leq k < 2^s$.
  \[
   {\it P_i \leq \left\{
                \begin{array}{ll}
                  \left(\frac{2k}{n}p_i + 1 - \frac{2k}{n} \right) p_i^s 
                
                  & \;\;\mbox{if } \;\;\;\;p_i > \frac{1}{2} \\
                  \\
                 \frac{1}{n}& \;\;\mbox{if }\;\;\;\;  p_i = \frac{1}{2}\\
                 \\
                \prod_{j=1}^{n-1} \left(\frac{2}{n-j+1}p_i + 1 - \frac{2}{n-j+1} \right)  & \;\;\mbox{if } \;\; \;\;p_i < \frac{1}{2}\\
                \end{array}
              \right.}
  \]
  
\end{theorem}
  \noindent {\bf Proof.} \hspace{1mm} The proof follows directly from Theorem 2 and the preceding Lemma. $\qquad \rule{2mm}{2mm}$ \\
  
  \begin{corollary}
 For a player with $p_i < \frac{1}{2}$, the probability of winning the tournament under any format is less than $\frac{1}{n}$.
  \end{corollary}

\end{document}